\documentstyle[%showkeys,
11pt,amsfonts]{article}

\newcommand {\rel} {{\mathbb R}}
\newcommand {\com} {{\mathbb C}}
\newcommand {\nat} {{\mathbb N}}

\newcommand {\ganz} {{\mathbb Z}}

\newtheorem{proposition}{Proposition}[section]
\newtheorem{theorem}{Theorem}[section]
\newtheorem{corollaryth}[theorem]{Corollary}
\newtheorem{lemma}[proposition]{Lemma}

\newcommand{\nc}{\newcommand}

\newcommand{\bcite}[1] {\cite{#1}}

\newcommand {\pr} {\bf}

\newcommand{\proof} {   \begin{flushright}
                        ///
                        \end{flushright}
                }

\newcommand{\partproof} {       \begin{flushright}
                                //
                                \end{flushright}
                }

\newcommand{\defin} { \hspace*{\fill} $\Box$ }

\newcommand{\foot}[1] { }
\renewcommand{\foot}[1] { \footnote{#1} }

\def    \mean   {{ \vec {\bf H} }}

\def	\d	{{ \ {\rm d} }}

\def    \W      {{ {\mathcal W} }}
\nc{\energ}[1]  {{ e_{#1} }}

\def	\minin	{{ {\mathcal M}_n }}
\newcommand{\miniz}[1] {{ {\mathcal M}_{#1} }}

\def	\metric	{{ {\cal M}et }}
\def	\modul	{{ {\cal M} }}
\newcommand{\coto}[1] {{ T^2_{#1} }}

\def	\poin	{{ poin }}

\def	\geu	{{ g_{euc} }}

\nc{\gpo}[1]	{{ g_{\poin}^{#1} }}
\nc{\gpu}[1]	{{ g_{\poin,{#1}} }}

\nc{\gpmo}[1]	{{ g_{\poin,m}^{#1} }}
\nc{\tgpmo}[1]	{{ \tilde g_{\poin,m}^{#1} }}

\nc{\gplo}[1]	{{ g_{\poin,\lambda}^{#1} }}

\def	\tp	{{ \pi }}

\def	\teich	{{ \cal T }}

\def    \diff	     {{ \stackrel{\approx}{\longrightarrow} }}

\nc{\doub}[1]{{ \ddot{#1} }}
\nc{\dd}{ \begin{displaymath} }
\nc{\df}{ \end{displaymath} }
\nc{\dcd}{ \begin{displaymath} \begin{array}{c}}
\nc{\dcf}{ \end{array} \end{displaymath} }
\nc{\ee}{ \begin{equation} }
\nc{\ef}{ \end{equation} }
\nc{\ad}{ \begin{array}{c} }
\nc{\af}{ \end{array} }

\begin{document}

%\begin{flushright}
%not for distribution
%\end{flushright}
\begin{center}
{\huge \bf Conformal Willmore Tori in $\rel^4$}
\\ \ \\
Tobias Lamm \\
Fakult\"at f\"ur Mathematik\\
Karlsruher Institut f\"ur Technologie (KIT)\\
Englerstra\ss e 2,
D-76131 Karlsruhe, Germany, \\
email: tobias.lamm@kit.edu \\
\ \\
Reiner M. Sch\"atzle \\
Fachbereich Mathematik der
Eberhard-Karls-Universit\"at T\"ubingen, \\
Auf der Morgenstelle 10,
%Geb\"aude C, 5 A 40,
D-72076 T\"ubingen, Germany, \\
email: schaetz@everest.mathematik.uni-tuebingen.de \\
 
\end{center}
\vspace{1cm}

\begin{quote}

{\bf Abstract:} For every two-dimensional torus $T^2$ and every $k\in \nat$, $k\ge 3$, we construct a conformal Willmore immersion $f:T^2\to \rel^4$ with exactly one point of density $k$ and Willmore energy $4\pi k$. Moreover, we show that the energy value $8\pi$ cannot be attained by such an immersion. Additionally, we characterize the branched double covers $T^2\to S^2 \times \{0\}$ as the only branched conformal immersions, up to M\"obius transformations of $\rel^4$, from a torus into $\rel^4$ with at least one branch point and Willmore energy $8\pi$. Using a perturbation argument in order to regularize a branched double cover, we finally show that the infimum of the Willmore energy in every conformal class of tori is less than or equal to $8\pi$.
\end{quote}

\vspace{1cm}

%\tableofcontents

%%%%%

\setcounter{equation}{0}

\section{Introduction} \label{intro}

The famous Willmore conjecture states that the minimum of the Willmore energy among all immersions $f:T^2 \to \rel^3$ from a two-dimensional torus is equal to $2\pi^2$ and it is attained by the stereographic image of the Clifford torus $S^1(\frac1{\sqrt{2}})\times S^1(\frac1{\sqrt{2}}) \subset S^3$. This conjecture was recently proved by Marques and Neves \cite{marquesneves}. We note that the Willmore energy $\W(f)$ for an immersion $f:\Sigma \to \rel^n$ of a Riemann surface $\Sigma$ is defined to be
\[
\W(f)=\frac14 \int \limits_{\Sigma} |H|^2 \d \mu,
\]
where $H$ and $\d \mu$ are the mean curvature resp. the induced area element of the immersion. Critical points of $\W$ are called Willmore immersions.

Another interesting question to ask is whether the infimum of the Willmore energy is also attained in every conformal class of tori, or even more generally, if the infimum is attained for every closed Riemann surface $\Sigma$ of genus $g\ge 1$. The immersions minimizing the Willmore energy in a fixed conformal class are called conformally constrained Willmore minimizers. 
The existence of these minimizers was established for some class of closed Riemann surfaces
in \bcite{kuw.schae.will7}.
This was extended in \bcite{kuw.li.mini} and \bcite{rivi.mini}
to Riemann surfaces $\ \Sigma\ $ which admit a conformal immersion
$\ f: \Sigma \rightarrow \rel^n\ $ with Willmore energy $\ \W(f) < 8 \pi\ $,
and in any case without energy restriction,
a branched conformally constrained minimizer was obtained.
Smoothness without energy restriction and in any codimension
was proved for unbranched conformally constrained minimizers
in \cite{kuw.schae.will7},
and for branched conformally constrained minimizers,
analyticity was proved in \cite{rivi.adv} under the assumption that either the genus of the surface is less than or equal to two or the Teichm\"uller class of the minimizing immersion is not hyperelliptic. Further, in \cite{rivi.adv} minimizing constrained by closed submanifolds of the Teichm\"uller space rather than a fixed class was considered.
Moreover, in two papers of Ndiaye and the second author \cite{ndiayeschaetzle1,ndiayeschaetzle2}, it was shown that the CMC-tori $S^1(r) \times S^1(\sqrt{1-r^2})\subset S^3$ minimize the Willmore energy in their own conformal class in arbitrary codimensions if $r\approx 1/\sqrt{2}$. 

In this paper we study tori in higher codimension $\ n \geq 4\ $.
Any torus is conformally equivalent
to a quotient $\ \coto{\omega} := \com / (\ganz + \omega \ganz)
\mbox{ with } \geu\ $ and
\begin{displaymath}
	\omega \in \modul := \{ a + ib\ |
	\ b > 0, 0 \leq a \leq 1/2,
	a^2 + b^2 \geq 1\ \},
\end{displaymath}
see \bcite{jo.rie} Theorem 2.7.1,
and we put
\begin{equation} \label{intro.mini}
	\minin(\omega)
	:= {\mathcal M}_{n,1}(\omega)
	:= \inf \{ \W(f)\ |\ f:
	\coto{\omega} \rightarrow \rel^n
	\mbox{ conformal } \}
\end{equation}
for $\ \omega \in \modul\ $. Our first main result is an existence statement for conformal Willmore tori in every conformal class with a prescribed energy value.
\\ \ \\
{\bf Theorem \ref{will.tori}}
{\it
%%%%%
For any conformal class $\ \omega \in \modul
\mbox{ and } k \in \nat_0, k \geq 3,\ $
there exist conformal Willmore immersions
$\ f_{\omega , k}: \coto{\omega} \rightarrow \rel^4\ $ with exactly one point of density $k$ and
\begin{displaymath}
	\W(f_{\omega , k}) = 4 k \pi
	\quad \mbox{for } k \geq 3.
\end{displaymath}%%%%%
}
\defin
\\ \ \\
We complement this result with a non-existence statement for conformal Willmore tori with at least one double point and Willmore energy $8\pi$.
\\ \ \\
{\bf Theorem \ref{nonexistence}}
{\it
%%%%%
For every torus $T^2$ there is no immersion $f_0:T^2\to \rel^4$ which has at least one double point and for which $\W(f_0)=8\pi$. %%%%%
}
\defin
\\ \ \\
A result similar to Theorem \ref{will.tori} for spheres was shown by Barbosa \cite{barbosa} and Calabi \cite{calabi}. For every integer $k\in \nat \backslash \{2\}$, they showed the existence of a superminimal immersion $\Phi:S^2 \to S^4$ with area $4k\pi$. Moreover, they showed that there is no minimal immersion with area $8\pi$. After stereographic projection at a point which is not in the image of the surfaces, one obtains Willmore spheres in $\rel^4$ with Willmore energy equal to $4k \pi$, $k\in \nat \backslash \{2\}$.
Note that on the other hand there exist Willmore spheres in $\rel^4$ or $S^4$ with Willmore energy $8\pi$ (for example the Whitney sphere), see e.g. \cite{montiel}. Furthermore, it was shown by Montiel \cite{montiel}, that the Willmore energy of every Willmore sphere $S^2 \rightarrow S^4$ (or $\rel^4$) has to be a multiple of $4\pi$, thereby generalizing the codimension one result of Bryant \cite{bryant84}. 

The particular implication of Theorem \ref{will.tori}
that for any conformal class $\ \omega \in \modul\ $
there exists a conformal Willmore immersion
$\ \coto{\omega} \rightarrow \rel^4\ $ is already known,
as Bryant showed in \bcite{brya.mini}
that any closed Riemann surface $\ \Sigma\ $ admits
a conformal minimal, even a superminimal, immersion $\ \Sigma \rightarrow S^4\ $,
which then is Willmore as well. Moreover, he constructed the immersions as a Twistor projection $T:\com P^3 \rightarrow S^4$ of a holomorphic horizontal curve $\Phi:\Sigma \to \com P^3$ and he showed that the Willmore energy of the superminimal immersion has to be a multiple of $4\pi$. Note however, that Bryant only obtains the existence of one such surface, whereas our result shows the existence of infinitely many Willmore immersions on every torus. 

Our construction of the conformal Willmore immersions is different to the one of Bryant, since we obtain our immersions via an inversion of suitable conformal minimal immersions in $\rel^4$ with ends of multiplicity one.
More precisely, we construct these immersions via a pair of meromorphic functions $(f,h):\coto{\omega} \rightarrow \rel^4$ with exactly $k\ge 3$ simple poles and no common branch points. The existence of these functions follows basically from the Riemann-Roch theorem. It then remains to show that by inverting the immersion $(f,h)$ one obtains an immersion as claimed in the theorem. 
In a remark after the proof of Theorem \ref{will.tori}, we show that all of our immersions are different from the one constructed by Bryant. 

A similar construction was employed previously by Weiner \cite{weiner}. He showed that for every compact Riemann surface $M$ of genus $p$ which is a holomorphic submanifold of $\com P^n$, for some $n$, of degree $d$, there exists a conformal immersion $f:M\to \rel^4$ with $\W(f)=4\pi d$ and which minimizes the Willmore energy in its regular homotopy class of immersions. Since every compact Riemann surface satisfies these assumptions for some value of $d$ and $n=3$ this implies the existence of a conformal Willmore immersion $f:M\to \rel^4$ for every such Riemann surface. We remark that the above immersions arise in the study of the equality case in an inequality derived by Wintgen \cite{wintgen}.

All constrained Willmore tori in $S^4$ were classified by Bohle \cite{bohle}. He showed that they are either superconformal or stereographic projections of minimal surfaces in $\rel^4$ with planar ends or the spectral curve has finite genus. The examples constructed by Bryant \cite{brya.mini} are superconformal and the minimal surfaces with flat ends constructed in Theorem \ref{will.tori} give rise to Willmore tori in $S^4$ via stereographic projection

We note that the situation is somehow different in codimension one. Namely, there do not exist minimal tori in $\rel^3$ with two or three embedded planar ends.
As finitely many embedded planar ends imply finite total curvature,
minimal tori in $\rel^3$ with two embedded planar ends
are excluded by \bcite{schoe.mini}.
For three embedded planar ends, this was ruled out by Kusner and Schmitt \cite{kusner.schmitt}. This shows in particular that Willmore tori in $\rel^3$ with Willmore energy $12\pi$ cannot be constructed via minimal surfaces.

In order to show the nonexistence result in Theorem \ref{nonexistence}, we invert the conformal immersion $f_0$ at the unit circle centered at one of its double points. It follows from a result of Weiner \cite{weiner} that the image $f$ is a conformal minimal immersion from a twice punctured torus and $\partial f$ is a meromorphic $\com^4$-valued one-form with poles of order two precisely at the punctures. Using the fact that every torus is conformally equivalent to $\coto{\omega}$ for some $\omega \in \modul$, we show that the doubly periodic $\partial_z f$ can be expressed as a linear combination of two suitably modified Weierstrass $\wp$-functions and a constant. Since $f$ is conformal, we are then able to show that we can reduce everything to codimension one, i.e. we construct out of $f$ a modified conformal minimal immersion $\tilde f$ from a twice punctured torus into $\rel^3$ with finite total curvature, which is a contradiction to the above mentioned result of Schoen \cite{schoe.mini}.
\\ \ \\
As an application and extension of the non-existence result for immersions from $T^2$ into $\rel^4$ with at least one double point and Willmore energy $8\pi$, we also classify all branched conformal immersions from $\coto{\omega}$, for every $\omega \in \modul$, into $\rel^4$ with at least one branch point and Willmore energy $8\pi$. We show that modulo M\"obius transformations these immersions are given by a branched double cover $\coto{\omega} \rightarrow S^2 \times \{0\}$. 

As already mentioned above, the existence of conformally constrained Willmore tori is known under the assumption $\minin(\omega)<8\pi$.  In our second main result we show that at least the non-strict inequality holds by perturbing a branched double cover. More precisely we have the following
\\ \ \\
{\bf Theorem \ref{infi.tori}}
{\it
%%%%%
For any conformal class $\ \omega \in \modul\ $, we have
\begin{displaymath}
	\miniz{4}(\omega) \leq 8 \pi,
\end{displaymath}
in particular $\ \minin \mbox{ is continuous for } n \geq 4\ $.%%%%%
}
\defin
\\ \ \\
We show this result by using a perturbation argument which slighly resembles the constructions of counterexamples to rigidity results  in \cite{lammschaetzle}. More precisely, close to the branch points of the branched double cover we add a small multiple of a suitably localized holomorphic function in the second component and the new immersion is conformal everywhere. The drawback of this construction is that this might change the conformal class of the torus and we cope with this problem by showing that the induced Teichm\"uller class is surjective if the perturbation is small enough. Altogether, this yields a sequence of conformal immersions from every torus whose Willmore energy converges to $8\pi$ from above.

Finally, we use the classification result for branched conformal immersions from $\coto{\omega}$ into $\rel^4$, in order to show that a similar construction cannot be done in codimension one.
As a byproduct we show that for every conformal class $\omega \in \modul$ with $\ \miniz{3}(\omega) \leq 8 \pi\ $, there exists a smooth, conformally constrained Willmore minimizer $f:\coto{\omega} \to \rel^3$, which improves the above mentioned existence results in \cite{kuw.li.mini} and \cite{rivi.mini} for tori in $\rel^3$. 
%%%%%
%%%%%

\setcounter{equation}{0}

\section{Willmore tori in higher codimension} \label{will}
In the first lemma of this section we show how one can construct a smooth immersion out of meromorphic functions with at most simple poles at the origin. A version of this result was already shown in \cite{weiner}.
\begin{lemma} \label{will.lem}

Let $\ f, f_1, f_2: B_1(0) - \{0\} \subseteq \com \rightarrow \com\ $
be smooth functions
which extend smoothly to the origin or
are holomorphic and have a simple pole at the origin,
and let $\ f_1\ $ be holomorphic with a simple pole at the origin.
Then
\begin{displaymath}
	\varphi := \frac{f}{|f_1|^2 + |f_2|^2}
\end{displaymath}
with $\ \varphi(0) := 0\ $ is smooth locally around the origin.
Moreover if $\ f\ $ is holomorphic with a simple pole at the origin,
then $\ D \varphi(0)\ $ has full rank.
\end{lemma}
{\pr Proof:} \\
Since $\ f_1\ $ has a pole at the origin,
we know $\ |f_1|^2 + |f_2|^2 \neq 0\ $
and $\ \varphi\ $ is well defined locally around $\ 0\ $.
We calculate for $\ z \neq 0\ $ close to $\ 0\ $
\begin{displaymath}
	\varphi(z) = \bar z
	\frac{z f(z)}{z f_1(z) \bar z \bar f_1(z) + z f_2(z) \bar z \bar f_2(z)}.
\end{displaymath}
Since $\ f, f_1, f_2\ $ are smooth or have at most a simple pole at $\ 0\ $,
the functions $\ h_i \mbox{ defined by } h_i(z) := z f_i(z), i = 0,1,2,\ $
are smooth.
We see
\begin{displaymath}
	\varphi(z) = \bar z \frac{h(z)}{|h_1(z)|^2 + |h_2(z)|^2}
\end{displaymath}
and, when observing that $\ h_1(0) \neq 0\ $,
as $\ f_1\ $ has a simple pole at $\ 0\ $,
we conclude that $\ \varphi\ $ is smooth at $\ 0\ $.

Further if $\ f\ $ is holomorphic with a simple pole at $\ 0\ $,
we know $\ h(0) \neq 0\ $
and calculate by standard Wirtinger calculus that
\begin{displaymath}
	\partial_z \varphi
	= \bar z \cdot \partial_z \Big( \frac{h}{|h_1|^2 + |h_2|^2} \Big),
\end{displaymath}
hence $\ \partial_z \varphi(0) = 0\ $,
and
\begin{displaymath}
	\partial_{\bar z} \varphi(0)
	= \frac{h(0)}{|h_1(0)|^2 + |h_2(0)|^2} \neq 0.
\end{displaymath}
Together
\begin{displaymath}
	\det D \varphi(0)
	= |\partial_z \varphi(0)|^2 - |\partial_{\bar z} \varphi(0)|^2 \neq 0,
\end{displaymath}
and $\ D \varphi(0)\ $ has full rank.
\proof
Next we show how the previous construction can be used in order to obtain smooth conformal Willmore immersions in $S^4$.
\begin{proposition} \label{will.prop}

Let $\ f_1, f_2: B_1(0) \rightarrow \com \cup \{ \infty \}\ $
be two meromorphic functions
with only simple poles and
\begin{equation} \label{will.prop.zero}
	[f_1' = 0] \cap [f_2' = 0] = \emptyset.
\end{equation}
Then for any M\"obius transformation $\ \Phi: \rel^4 \cup \{ \infty \} \diff S^4\ $, i.e. a composition of a M\"obius transformation $\ \varphi:S^4 \to S^4\ $ and the inverse of an arbitrary stereographic projection $\ T:S^4 \diff \rel^4 \cup \{ \infty \} \ $,
the map
\begin{displaymath}
	\Phi \circ (f_1,f_2): B_1(0) \rightarrow S^4
\end{displaymath}
is a smooth conformal Willmore immersion.
\end{proposition}
{\pr Proof:} \\
If both $\ f_1, f_2\ $ are holomorphic,
then $\ (f_1,f_2): B_1(0) \rightarrow \com^2\ $ is smooth as well
and with pull-back metric
\begin{displaymath}
	g := (f_1,f_2)^* \geu = f_1^* \geu + f_2^* \geu
	= (|f_1'|^2 + |f_2'|^2) \geu,
\end{displaymath}
as $\ f_1, f_2\ $ are conformal by holomorphy.
Then by (\ref{will.prop.zero}),
we see that $\ (f_1,f_2)\ $ is a conformal immersion.

Moreover by holomorphy, $\ f_1, f_2\ $ are harmonic,
hence by conformal invariance of the Laplacian
\begin{displaymath}
	\mean_{(f_1,f_2)} = \Delta_g (f_1,f_2) = 0,
\end{displaymath}
and we conclude that $\ (f_1,f_2)\ $ is minimal,
in particular a Willmore immersion.
By conformal invariance $\ \Phi \circ (f_1,f_2)\ $ is a Willmore immersion as well,
hence proving the proposition in the absence of poles.

Now the poles of $\ f_1, f_2\ $ do not accumulate in $\ B_1(0)\ $
by the definition of meromorphic functions,
and we assume that $\ f_1 \mbox{ or } f_2\ $
has a simple pole at $\ 0\ $.
Instead of considering a M\"obius transformation $\ \Phi: \rel^4 \cup \{ \infty \} \diff S^4\ $,
we consider the inversion $\ I \mbox{ of } \rel^4 - \{0\}
\mbox{ given by } (z,w) \mapsto (z,w) / (|z|^2 + |w|^2)\ $.
By the previous lemma, we already know that
\begin{displaymath}
	I \circ (f_1,f_2)
	= \frac{(f_1,f_2)}{|f_1|^2 + |f_2|^2}
\end{displaymath}
is smooth locally around $\ 0\ $
and that $\ D(I \circ (f_1,f_2))(0)\ $ has full rank,
hence $\ I \circ (f_1,f_2)\ $ is a smooth immersion locally around $\ 0\ $.
Next $\ I \circ (f_1,f_2)\ $ is
a smooth Willmore immersion in a punctured disc of $\ 0\ $
by the argument above,
hence $\ I \circ (f_1,f_2)\ $ is a smooth Willmore immersion at $\ 0\ $ as well.
Since every M\"obius transformation $\ \Phi: \rel^4 \cup \{ \infty \} \diff S^4\ $ can be written as $\Phi' \circ I$, where $\ \Phi': \rel^4 \cup \{ \infty \} \diff S^4\ $ is again a M\"obius transformation, this concludes the proof of the proposition.
\proof
Combining the above two results with the Riemann-Roch theorem, which yields the existence of a meromorphic function with $m\ge 2$ simple poles on every torus, we are now able to show the first main theorem.
\begin{theorem} \label{will.tori}

%%%%%
For any conformal class $\ \omega \in \modul
\mbox{ and } k \in \nat_0, k \geq 3,\ $
there exist conformal Willmore immersions
$\ f_{\omega , k}: \coto{\omega} \rightarrow \rel^4\ $ with exactly one point of density $k$ and
\begin{displaymath}
	\W(f_{\omega , k}) = 4 k \pi.
\end{displaymath}%%%%%
\end{theorem}
{\pr Proof:} \\
By the Riemann-Roch theorem,
there exists a meromorphic function $\ f: \coto{\omega} \rightarrow S^2\ $
with $\ m\ $ simple poles for $\ m \geq 2\ $,
%say at $\ p_1 \neq \ldots \neq p_m \in \coto{\omega}\ $,
see \bcite{jo.rie} Theorem 5.4.1,
in particular $\ f\ $ is of degree $\ m\ $.
Clearly this is a branched conformal immersion.
Considering a meromorphic function $\ h\ $ with two simple poles 
and the covering projection $\ \pi: \com \rightarrow \coto{\omega}\ $,
any translation in $\ \com \mbox{ by } v \in \com\ $
induces a conformal automorphism $\ \tau_v \mbox{ of } \coto{\omega}\ $.
Since $\ f \mbox{ and } h\ $ have only finitely many branch points,
we may replace $\ h \mbox{ by } h \circ \tau_v
\mbox{ for appropriate } v \in \com\ $ such that
\begin{displaymath}
	[d f = 0] \cap [d h = 0] = \emptyset.
\end{displaymath}
Moreover we may assume that $\ h\ $
has no branching at the poles of $\ f\ $
and that none of the poles of $\ f\ $
coincides with a pole of $\ h\ $.
In this case
\begin{displaymath}
	k := \# \{ \mbox{ poles of } f \mbox{ or } h\ \} = m + 2.
\end{displaymath}
As $\ m \geq 2\ $ was arbitrary,
we can achieve for $\ k\ $ any integer $\ \geq 4\ $.
To achieve also $\ k = 3\ $,
we consider $\ m = 2\ $
and replace $\ h \mbox{ by } (h - h(p))^{-1}
\mbox{ for some pole } p \mbox{ of } f\ $,
hence $\ f \mbox{ and } h\ $ have at least one pole in common
and the poles of $\ h\ $ are still simple,
as $\ p\ $ is not a branch point of $\ h\ $ by above.
Next, if the two poles of $f$ and $h$ were the same,
then $h = \alpha f + \beta$ by the Riemann-Roch theorem,
which is impossible, since we already know
that $f$ and $h$ do not have common branch points.
Hence we get $\ k = 3\ $.

By the previous proposition,
$\ \Phi \circ (f,h): \coto{\omega} \rightarrow S^4\ $
is a smooth conformal Willmore immersion
for any M\"obius transformation $\ \Phi: \rel^4 \cup \{ \infty \} \diff S^4\ $.
Moreover since $\ f,h\ $ are holomorphic
in $\ \coto{\omega} - \{ \mbox{ poles of } f \mbox{ or } h\ \}\ $,
hence harmonic, we see that $\ (f,h): \coto{\omega} - \{ \mbox{ poles of } f \mbox{ or } h\ \}
\rightarrow \com^2 \cong \rel^4\ $ is minimal, hence
\begin{displaymath}
	\W( (f,h) ) = 0.
\end{displaymath}
Clearly by above, the preimages of infinity under $\ (f,h)\ $
are exactly the poles of $\ f \mbox{ and } h\ $,
and since for the inversion $\ I \mbox{ of } \rel^4 \cong \com^2\ $
the map $\ I \circ (f,h)\ $ is a smooth immersion
locally around the poles of $\ f \mbox{ and } h\ $ by the previous lemma,
we see by a standard calculation and the Gau\ss-Bonnet theorem,
%see also \bcite{lamm.nguy.conf} Theorem 2.3,
that
\begin{displaymath}
	\W(I \circ (f,h)) =  \W( (f,h) )
	+ 4 \pi \cdot \#( (f,h)^{-1}(\infty) )
	= 4 k \pi,
\end{displaymath}
where we also assumed without loss of generality that $0\notin (f,h)(\coto\omega)$. 
Letting $f_{\omega,k}:=I\circ (f,h)$ concludes the proof of the proposition.
\proof
{\large \bf Remark:} \\ 
It follows from \cite{Friedrich} that the immersion $(f,h):\coto{\omega} - \{ \mbox{ poles of } f \mbox{ or } h\ \}\  \rightarrow \com^2 \cong \rel^4$ constructed in Theorem \ref{will.tori} is superminimal. From the above results we get that $\Phi_N^{-1} \circ (f,h) =\Phi_S^{-1}\circ I \circ (f,h):\coto{\omega} \to S^4$ is a smooth conformal Willmore immersion. Here $\Phi_N$ resp. $\Phi_S$ denote the stereographic projections from the north resp. south pole.
Now we claim that the immersion $\Phi_N^{-1} \circ (f,h)$ is not minimal:

If we assume that both $(f,h)$ and $\Phi_N^{-1} \circ (f,h)$ are minimal, then it follows from the transformation formula for the mean curvature under conformal changes
\begin{equation}\label{trafoH}
\hat H=\lambda^{-2} (H- \lambda^{-1} (\nabla \lambda)^\perp),
\end{equation}
where $\hat H$ is the mean curvature of the immersion with target $(\rel^4,\lambda^2 \geu)$, $\lambda= 2/(1+|x|^2)$ and $v^\perp$ denotes the normal component of a vector $v\in \rel^4$, that $(f,h)$ is a minimal cone. Since the intersection of every 2-dimensional minimal cone with the unit sphere $S^3$ of $\rel^4$ is a geodesic, it follows that $(f,h)$ has to be a plane and $\Phi_N^{-1} \circ (f,h)\subset S^4$ a geodesic sphere, which contradicts the fact that our surfaces are tori.

Finally, we note that every M\"obius transformation $\Phi: \rel^4 \cup \{ \infty \} \diff S^4$ can be written as $\Phi = R \circ \Phi_N^{-1} \circ T$, where $R$ is a rotation of $S^4$ und $T$ is a composition of rotations, dilations and translations of $\rel^4$. This follows from choosing a rotation $R$ of $S^4$ so that $R(N)=\Phi(\infty)$. The map $T := \Phi_N \circ R^{-1} \circ \Phi :\rel^4 \cup \{ \infty \} \to \rel^4 \cup \{ \infty \}$ is then a M\"obius transformation with $T(\infty)=\infty$. 
Now the rotation $R$ maps minimal surfaces onto minimal surfaces and great circles onto great circles and $T$ also maps minimal surfaces onto minimal surfaces. 

Hence the above argument shows that $\Phi \circ (f,h)$ cannot be minimal in $S^4$ for every M\"obius transformation $\Phi: \rel^4 \cup \{ \infty \} \diff S^4$. In particular, all these surfaces are genuinely different from the ones constructed by Bryant in \cite{brya.mini}.
\defin
\\ \\
%%%%%
The Willmore immersions constructed in Theorem \ref{will.tori}
with Willmore energy $\ 4 \pi k, k \geq 3\ $,
all have a point of density $\ k\ $.
We show that this is impossible for $\ k = 2\ $. 
\begin{theorem}\label{nonexistence}
For every torus $T^2$ there is no immersion $f_0:T^2\to \rel^4$ which has at least one double point and for which $\W(f_0)=8\pi$. 
\end{theorem}
{\pr Proof:} \\
Indeed if there is an immersion $\ f_0: T^2 \rightarrow \rel^4\ $
with $\ \W(f_0) = 8 \pi\ $ and at least one double point,
say $\ f_0(p_1) = f_0(p_2) = 0 \mbox{ for some } p_1 \neq p_2 \in T^2\ $
after translation,
then by \bcite{weiner} Proposition 2
its inversion $\ f := f_0 / |f_0|^2: T^2 - \{ p_1 , p_2 \} \rightarrow \rel^4\ $
is a minimal immersion
and $\ \partial f\ $ is a meromorphic $\ \com^4-\mbox{ valued } 1-$form on $\ T^2\ $
with poles precisely at $\ p_1 \neq p_2\ $
and each pole is of order two.
By conformal equivalence $\ (T^2,f_0^* \geu) \cong \coto{\omega}
\mbox{ for some } \omega \not\in \rel\ $,
we may consider $\ f_0 \mbox{ and } f\ $
as conformal immersions doubly periodic with respect to $\ \Gamma := \ganz + \omega \ganz\ $
on $\ \com \mbox{ respectively on } \com - \{ p_1 + \Gamma \neq p_2 + \Gamma \}, p_1, p_2 \in \com\ $,
and $\ \partial_z f = (\partial_1 f - i \partial_2 f) / 2\ $
is a doubly periodic meromorphic $\ \com^4-$valued function
with poles precisely at $\ p_1 + \Gamma \neq p_2 + \Gamma\ $
and each pole is of order two.

Next the Weierstrass $\ \wp-$function on $\ \com / \Gamma\ $ given by
\begin{displaymath}
	\wp(z) := \frac{1}{z^2} + \sum \limits_{\gamma \in \ganz + \omega \ganz - \{0\}}
	\Big( \frac{1}{(z-\gamma)^2} - \frac{1}{\gamma^2} \Big),
\end{displaymath}
see \bcite{ahl} \S 7.3.1, has a pole of order two
with residue zero precisely at $\ \Gamma\ $.
Then
\begin{equation} \label{will-2.weier}
	\wp_l := \wp(. - p_l) - \wp(p_{3-l} - p_l),
	\quad l = 1,2,
\end{equation}
has a pole of order two
with residue zero precisely at $\ p_l + \Gamma\ $
and vanishes at $\ p_{3-l} + \Gamma\ $.
Clearly as $\ \wp\ $ is an even function by above,
\begin{equation} \label{will-2.weier-const}
	\wp(p_2-p_1) = \wp(p_1 - p_2).
\end{equation}
Then
\begin{equation} \label{will-2.weier-prod}
	\wp_1 \wp_2 \mbox{ has at most simple poles at } p_1 + \Gamma, p_2 + \Gamma,
\end{equation}
and for appropriate $\ a,b \in \com^4\ $,
also $\ \partial_z f - a \wp_1 - b \wp_2\ $
has at most simple poles at $\ p_1 + \Gamma, p_2 + \Gamma\ $.
Clearly as $\ \partial_z f\ $ has a pole of order two
at both $\ p_1 + \Gamma, p_2 + \Gamma\ $,
we have
\begin{equation} \label{will-2.ab-null}
	a,b \neq 0.
\end{equation}
As in the proof of Proposition \ref{will.tori}
by the Riemann-Roch theorem,
there exists a doubly periodic meromorphic function $\ w: \com \rightarrow S^2\ $
with simple poles precisely at $\ p_1 + \Gamma \neq p_2 + \Gamma\ $,
see \bcite{jo.rie} Theorem 5.4.1.
We see by periodicity for the fundamental domain
\begin{displaymath}
	I := \{ s + t \omega\ |\ s,t \in [0,1]\ \}
\end{displaymath}
and $\ \xi \in \com \mbox{ with } p_l \not\in \Gamma + \partial(\xi + I)\ $ that
\begin{equation} \label{will-2.residue}
	0 = \frac{1}{2 \pi i} \int \limits_{\partial(\xi + I)} w(\zeta) \d \zeta
	= Res(w,p_1) + Res(w,p_2),
\end{equation}
and we may assume that $\ Res(w,p_1) = 1, Res(w,p_2) = -1\ $.
Likewise $\ Res(\partial_z f,p_1) + Res(\partial_z f,p_2) = 0\ $,
and for $\ c \in \com^4\ $
we get that $\ \partial_z f - a \wp_1 - b \wp_2 - c w\ $
is doubly periodic meromorphic without poles, hence is constant and
\begin{equation} \label{will-2.fz}
	\partial_z f = a \wp_1 + b \wp_2 + c w + d
\end{equation}
for some $\ d \in \com^4\ $.
As $\ \partial_z f\ $ is the derivative
of the real function $\ f \mbox{ in } B_\varrho(p_l) - \{p_l\}
\mbox{ for small } \varrho > 0\ $,
the period of $\ \partial_z f \mbox{ around } p_l\ $
has to be purely imaginary, that is
\begin{displaymath}
	i \rel^4 \ni \int \limits_{\partial B_\varrho(p_l)} \partial_z f(\zeta) \d \zeta
	= \int \limits_{\partial B_\varrho(p_l)} (a \wp_1 + b \wp_2 + c w + d)(\zeta) \d \zeta =
\end{displaymath}
\begin{displaymath}
	= 2 \pi i c Res(w,p_l) = 2 \pi i (-1)^{l-1} c,
\end{displaymath}
when recalling that the residues of $\ \wp_1, \wp_2 \mbox{ and } d\ $ vanish
and that $\ Res(w,p_l) = (-1)^{l-1}\ $.
We conclude that $\ c \in \rel^4\ $.

Introducing $\ \langle z , w \rangle = \sum_{j=1}^4 z_j w_j
\mbox{ for } z,w \in \com^4\ $,
conformality of $\ f\ $ reads as
\begin{displaymath}
	0 = \langle \partial_z f , \partial_z f \rangle =
\end{displaymath}
\begin{displaymath}
	= \langle a,a \rangle \wp_1^2 + \langle b,b \rangle \wp_2^2 +
\end{displaymath}
\begin{displaymath}
	+ 2 \langle a,c \rangle \wp_1 w + 2 \langle b,c \rangle \wp_2 w +
\end{displaymath}
\begin{displaymath}
	+ 2 \langle a,d \rangle \wp_1 + 2 \langle b,d \rangle \wp_2
	+ \parallel c \parallel^2 w^2 +
\end{displaymath}
\begin{equation} \label{will-2.conf}
	+ 2 \langle a,b \rangle \wp_1 \wp_2
	+ 2 \langle c,d \rangle w
	+ \langle d,d \rangle.
\end{equation}
As $\ \wp_l^2\ $ is the only function with a pole of order four at $\ p_l\ $,
we get
\begin{equation} \label{will-2.ab}
	\langle a,a \rangle, \langle b,b \rangle = 0.
\end{equation}
Then $\ \wp_l w\ $ is the only apprearing function with a pole of order three at $\ p_l\ $,
and we get
\begin{equation} \label{will-2.ab-c}
	\langle a,c \rangle, \langle b,c \rangle = 0.
\end{equation}
Next by (\ref{will-2.weier-prod}),
the only remaining functions with poles of order two
are $\ \wp_1, \wp_2 \mbox{ and } w^2\ $.
$\ \wp_l\ $ has one pole of order two precisely at $\ p_l\ $
and with leading term $\ 1 / (z-p_l)^2\ $.
Since $\ w\ $ has two simple poles at $\ p_1, p_2\ $
with residue $\ \pm 1\ $,
we see that $\ w^2\ $ has two poles of order two at both $\ p_1, p_2\ $
and with leading term $\ 1 / (z-p_l)^2\ $.
Therefore
\begin{equation} \label{will-2.ad}
\begin{array}{c}
	2 \langle a,d \rangle + \parallel c \parallel^2 = 0, \\
	2 \langle b,d \rangle + \parallel c \parallel^2 = 0. \\
\end{array}
\end{equation}
We claim
\begin{equation} \label{will-2.c}
	c = 0.
\end{equation}
Indeed if $\ c \neq 0\ $,
we may assume after a rotation and a homothety of $\ \rel^4\ $
that $\ c = e_4\ $, as $\ c\ $ is real.
Then $\ a_4 = b_4 = 0\ $ by (\ref{will-2.ab-c}), and by (\ref{will-2.fz})
\begin{displaymath}
	\partial_z f_4 = w + d_4
\end{displaymath}
has simple poles at $\ p_1 + \Gamma, p_2 + \Gamma\ $ and expands to
\begin{displaymath}
	\partial_z f_4(p_1 + z) = \frac{1}{z} + \varphi(z)
	\quad \mbox{locally around } p_1
\end{displaymath}
with $\ \varphi\ $ holomorphic.
Choosing a holomorphic $\ \psi \mbox{ with } \psi' = \varphi
\mbox{ locally around } p_1\ $, we calculate
\begin{displaymath}
	2 \partial_z( \log |z| + Re(\psi)(z))
	= \frac{1}{z} + \psi'(z) = \partial_z f_4
	\quad \mbox{in } B_\varrho(0) - \{0\},
\end{displaymath}
hence $\ \partial_z(f_4(p_1 + .) - \log |.| - Re(\psi)) = 0\ $ and
\begin{displaymath}
	f_4(p_1 + z) = \log |z| + Re(\psi(z))
	\quad \mbox{in } B_\varrho(0) - \{0\}
\end{displaymath}
for appropriate holomorphic $\ \psi\ $.
We conclude
\begin{displaymath}
	f_4(p_1 + z) \asymp \log |z|
	\quad \mbox{for } z \rightarrow 0.	
\end{displaymath}
Since $\ f_0 = f / |f|^2\ $ is a smooth immersion vanishing at $\ p_1\ $,
we get $\ |f_0(p_1 + z)| \asymp |z|\ $ and
\begin{displaymath}
	|f(p_1+z)| \asymp 1/|z|
	\quad \mbox{for } z \rightarrow 0.	
\end{displaymath}
Together
\begin{displaymath}
	|f_{0,4}(p_1+z)| = |f_4(p_1+z)| / |f(p_1+z)|^2
	\asymp |z|^2 \log |z|
	\quad \mbox{for } z \rightarrow 0.
\end{displaymath}
Then by smoothness of $\ f_{0,4}\ $,
we get $\ f_{0,4}(p_1) = 0, D f_{0,4}(p_1) = 0\ $,
hence $\ f_{0,4}(p_1+z) \leq C |z|^2\ $.
This is a contradiction,
and (\ref{will-2.c}) is proved.
\partproof
Combining (\ref{will-2.conf}), (\ref{will-2.ab}), (\ref{will-2.ad})
and (\ref{will-2.c}) yields
\begin{equation} \label{will-2.abd}
	2 \langle a,b \rangle \wp_1 \wp_2 + \langle d,d \rangle = 0.
\end{equation}
Before proceeding, we prove that
\begin{equation} \label{will-2.d-null}
	a_j + b_j = 0 \Longrightarrow d_j = 0
	\quad \mbox{for any } j = 1,2,3,4.
\end{equation}
Indeed considering $\ j = 4\ $
and as $\ \partial_z f_4\ $ is the derivative
of a doubly periodic real function,
the period of $\ \partial_z f_4\ $ with respect to any closed path
in $\ \coto{\omega} - \{ p_1 , p_2 \}\ $
has to be purely imaginary.
This reads for the paths $\ [\xi,\xi + \omega_k],
\omega_1 = 1, \omega_2 = \omega\ $,
which are closed in $\ \coto{\omega}\ $,
and appropriate $\ \xi \in \com\ $ with
\begin{displaymath}
	p_l \not\in \Gamma + \partial(\xi - p_1 + J)
\end{displaymath}
where
\begin{displaymath}
	J := \{ s \omega_k + t(p_1 - p_2)\ |\ s,t \in [0,1]\ \} 
\end{displaymath}
see below,
when using (\ref{will-2.fz}) and (\ref{will-2.c}) that
\begin{equation} \label{will-2.aux}
	i \rel \ni \int \limits_{[\xi,\xi + \omega_k]} \partial_z f_4(\zeta) \d \zeta
	= a_4 \int \limits_{[\xi,\xi + \omega_k]} \wp_1(\zeta) \d \zeta
	+ b_4 \int \limits_{[\xi,\xi + \omega_k]} \wp_2(\zeta) \d \zeta
	+ d_4 \omega_k.
\end{equation}
By definition of $\ \wp_l\ $ in (\ref{will-2.weier})
and by (\ref{will-2.weier-const}), we see
\begin{displaymath}
	\int \limits_{[\xi,\xi + \omega_k]} (\wp_1 - \wp_2)(\zeta) \d \zeta
	= \int \limits_{[\xi,\xi + \omega_k]} (\wp(\zeta - p_1) - \wp(\zeta - p_2)) \d \zeta =
\end{displaymath}
\begin{displaymath}
	= \int \limits_{[\xi-p_1,\xi-p_1+\omega_k] - [\xi-p_2,\xi-p_2+\omega_k]} \wp(\zeta) \d \zeta =
\end{displaymath}
\begin{displaymath}
	= \int \limits_{[\xi-p_1,\xi-p_1+\omega_k] + [\xi-p_1+\omega_k,\xi-p_2+\omega_k]
	+ [\xi-p_2+\omega_k,\xi-p_2] + [\xi-p_2,\xi-p_1]} \wp(\zeta) \d \zeta = 0
\end{displaymath}
\begin{equation} \label{will-2.weier-residue}
	= \int \limits_{\partial(\xi-p_1 + J)} \wp(\zeta) \d \zeta = 0
\end{equation}
by periodicity and meromorphy of $\ \wp\ $
and recalling that all residues of $\ \wp\ $ vanish.
Putting
\begin{displaymath}
	\sigma_k :=  \int \limits_{[\xi,\xi + \omega_k]} \wp_l(\zeta) \d \zeta
\end{displaymath}
independent of $\ l = 1,2\ $, we rewrite (\ref{will-2.aux}) into
\begin{displaymath}
	(a_4 + b_4) \sigma_k + d_4 \omega_k \in i \rel.
\end{displaymath}
Now if $\ a_4 + b_4 = 0\ $
then $\ d_4 \omega_k \in i \rel\ $,
hence $\ d_4 = 0\ $,
as $\ \omega_2 / \omega_1 = \omega \not\in \rel\ $,
and (\ref{will-2.d-null}) is proved.
\partproof
We continue with the case
\begin{equation} \label{will-2.ab-0}
	\langle a,b \rangle = 0.
\end{equation}
From (\ref{will-2.ab}),
we see that $\ Re(a) \perp Im(a), |Re(a)| = |Im(a)|\ $.
As $\ a \neq 0\ $ by (\ref{will-2.ab-null}),
we may assume after a rotation and a homothety of $\ \rel^4\ $
that $\ a = e_1 - i e_2\ $.
Then by (\ref{will-2.ab-0})
\begin{displaymath}
	0 = \langle a,b \rangle
	= b_1 - i b_2,
\end{displaymath}
and by (\ref{will-2.ab})
\begin{displaymath}
	0 = \langle b,b \rangle
	= b_1^2 + b_2^2 + b_3^2 + b_4^2
	= b_3^2 + b_4^2,
\end{displaymath}
hence likewise $\ b_4 = -i b_3 \mbox{ or } b_3 - i b_4 = 0\ $
after reflection at $\ \rel^3 \times \{0\}\ $.
Next combining (\ref{will-2.ad}) and (\ref{will-2.c}),
we get $\ \langle a,d \rangle, \langle b,d \rangle = 0\ $.
Firstly
\begin{displaymath}
	0 = \langle a,d \rangle = d_1 - i d_2,
\end{displaymath}
and secondly by above
\begin{displaymath}
	0 = \langle b,d \rangle
	= b_1 d_1 + b_2 d_2 + b_3 d_3 + b_4 d_4 =
\end{displaymath}
\begin{displaymath}
	= b_1 d_1 + (-i) b_1 (-i) d_1 + b_3 d_3 - i b_3 d_4
	= b_3 (d_3 - i d_4),
\end{displaymath}
hence
\begin{displaymath}
	d_3 - i d_4 = 0,
\end{displaymath}
if $\ b_3 \neq 0\ $.
Otherwise if $\ b_3 = 0\ $ then also $\ b_4 = 0\ $,
hence $\ a_j + b_j = 0 \mbox{ for } j = 3,4\ $
and $\ d_j = 0\ $ by (\ref{will-2.d-null}).
In both cases we have $\ d_3 - i d_4 = 0\ $.

Combining the above, we see
\begin{displaymath}
	\partial_z f_1 - i \partial_z f_2, \partial_z f_3 - i \partial_z f_4 = 0,
\end{displaymath}
hence, as $\ f\ $ is real,
\begin{displaymath}
	0 = \overline{\partial_z f_1 - i \partial_z f_2}
	= \partial_{\bar z}(f_1 + i f_2)
\end{displaymath}
and $\ h_1 := f_1 + i f_2 \mbox{ and likewise } h_2 := f_3 + i f_4\ $
are doubly periodic meromorphic functions
with poles only at $\ p_1, p_2\ $.
Clearly by elementary function theory
\begin{displaymath}
	h_1' = 2 \partial_z Re(h_1) = 2 \partial_z f_1
	= 2i \partial_z Im(h_1) = 2i \partial_z f_2
\end{displaymath}
and likewise for $\ h_2 $, hence
\begin{displaymath}
	\partial_z f = \frac{1}{2} (h_1' , -i h_1' , h_2' , -i h_2').
\end{displaymath}
Therefore $\ h_1', h_2'\ $ have poles of order at most two at $\ p_1, p_2\ $,
and $\ h_1, h_2\ $ have at most simple poles at $\ p_1, p_2\ $.
As (\ref{will-2.residue}) by periodicity,
the two residues of $\ h_l\ $ add up to zero, hence
\begin{displaymath}
	h_l = \alpha_l w + \beta_l
\end{displaymath}
for appropriate $\ \alpha_l, \beta_l \in \com, l = 1,2\ $.
Clearly $\ w\ $ has branch points, that is $\ w'(p) = 0\ $
for at least one $\ p \in \com - \{ p_1 + \Gamma, p_2 + \Gamma \}\ $.
Then $\ h_1'(p), h_2'(p) = 0\ $
and by above $\ \partial_z f(p) = 0\ $.
This contradicts our assumption
that $\ f_0\ $ is an immersion on $\ \com\ $
and $\ f\ $ is an immersion on $\ \com - \{ p_1 + \Gamma, p_2 + \Gamma \}\ $.
Therefore the case $\ \langle a,b \rangle = 0\ $ is impossible.
\partproof
In the remaining case when
\begin{equation} \label{will-2.ab-not0}
	\langle a,b \rangle \neq 0,
\end{equation}
we reduce to codimension one.
After a rotation of $\ \rel^4\ $, we may assume that
\begin{displaymath}
	Re(a+b), Im(a+b) \in \rel^2 \times \{0\},
\end{displaymath}
hence
\begin{equation} \label{will-2.ab-add}
	a_3 + b_3, a_4 + b_4 = 0
\end{equation}
and further $\ d_3, d_4 = 0\ $ by (\ref{will-2.d-null}).

We define
\begin{displaymath}
	\tilde f := (f_1,f_2,\tilde f_3),
	\tilde a := (a_1,a_2,\tilde a_3),
	\tilde b := (b_1,b_2,\tilde b_3),
	\tilde d := (d_1,d_2,0)
\end{displaymath}
for appropriate $\ \tilde f_3, \tilde a_3, \tilde b_3\ $ to be chosen below.
Clearly by (\ref{will-2.ad}) and (\ref{will-2.c})
\begin{displaymath}
	\langle \tilde a,\tilde d \rangle = \langle a,d \rangle = 0,
	\quad
	\langle \tilde b,\tilde d \rangle = \langle b,d \rangle = 0.
\end{displaymath}
Choosing $\ \tilde a_3 \in \com\ $ with
\begin{displaymath}
	\tilde a_3^2 = a_3^2 + a_4^2
\end{displaymath}
and $\ \tilde b_3 = -\tilde a_3\ $,
hence by (\ref{will-2.ab-add})
\begin{displaymath}
	\tilde b_3^2 = \tilde a_3^2
	= a_3^2 + a_4^2 = b_3^2 + b_4^2,
\end{displaymath}
we see $\ \tilde a_3 + \tilde b_3 = 0\ $ and by (\ref{will-2.ab})
\begin{displaymath}
	\langle \tilde a , \tilde a \rangle
	= a_1^2 + a_2^2 + \tilde a_3^2
	= a_1^2 + a_2^2 + a_3^2 + a_4^2
	= \langle a,a \rangle = 0
\end{displaymath}
and likewise $\ \langle \tilde b , \tilde b \rangle = 0\ $.
Further using (\ref{will-2.ab-add})
\begin{displaymath}
	\langle \tilde a , \tilde b \rangle
	= a_1 b_1 + a_2 b_2 + \tilde a_3 \tilde b_3
	= a_1 b_1 + a_2 b_2 - \tilde a_3^2 =
\end{displaymath}
\begin{equation} \label{will-2.ab-tilde}
	= a_1 b_1 + a_2 b_2 - a_3^2 - a_4^2
	= a_1 b_1 + a_2 b_2 + a_3 b_3 + a_4 b_4
	= \langle a,b \rangle.
\end{equation}
Also we want to establish that
\begin{equation} \label{will-2.tilde-not0}
	\tilde a_3 \neq 0.
\end{equation}
Indeed if $\ \tilde a_3 = 0\ $
then $\ a_3^2 + a_4^2 = \tilde a_3^2 = 0\ $
and by (\ref{will-2.ab})
\begin{displaymath}
	a_1^2 + a_2^2 = -a_3^2 - a_4^2 = 0,
\end{displaymath}
hence $\ a_2 = \pm i a_1\ $.
Next by (\ref{will-2.ad}) and (\ref{will-2.c})
\begin{displaymath}
	0 = \langle a,d \rangle
	= a_1 d_1 + a_2 d_2
	= a_1 (d_1 \pm i d_2).
\end{displaymath}
If $\ a_1 = 0\ $ then also $\ a_2 = \pm i a_1 = 0\ $ and by (\ref{will-2.ab-add})
\begin{displaymath}
	\langle a,b \rangle
	= a_3 b_3 + a_4 b_4 = -a_3^2 - a_4^2 = 0,
\end{displaymath}
contrary to our assumption (\ref{will-2.ab-not0}).
Therefore $\ a_1 \neq 0\ $,
and we get $\ d_1 \pm i d_2 = 0\ $.
This yields
\begin{displaymath}
	\langle d,d \rangle = d_1^2 + d_2^2 = 0
\end{displaymath}
and by (\ref{will-2.abd}) that
\begin{displaymath}
	\langle a,b \rangle \wp_1 \wp_2 = 0.
\end{displaymath}
As $\ \wp_1, \wp_2 \not\equiv 0\ $ are meromorphic functions,
we get $\ \langle a,b \rangle = 0\ $,
contrary to our assumption (\ref{will-2.ab-not0}),
and (\ref{will-2.tilde-not0}) is proved.

We want to define $\ \tilde f_3\ $ in such a way that
\begin{equation} \label{will-2.f3}
	\partial_z \tilde f_3 = \tilde a_3 \wp_1 + \tilde b_3  \wp_2
	= \tilde a_3 (\wp_1 - \wp_2).
\end{equation}
Then
\begin{displaymath}
	\partial_z \tilde f = \tilde a \wp_1 + \tilde b \wp_2 + \tilde d
\end{displaymath}
and as in (\ref{will-2.conf})
\begin{displaymath}
	\langle \partial_z \tilde f , \partial_z \tilde f \rangle
	= 0,
\end{displaymath}
that is $\ \tilde f: \com - \{ p_1 + \Gamma , p_2 + \Gamma \} \rightarrow \rel^3\ $
is weakly conformal.

We turn to (\ref{will-2.f3}).
Now $\ \wp_1 - \wp_2\ $ is a doubly periodic meromorphic function
with poles precisely at $\ p_1 + \Gamma, p_2 + \Gamma\ $
and these poles are of order two with residue zero,
more precisely with princial part
\begin{displaymath}
	(-1)^{l-1} (z-p_l)^{-2}
	\quad \mbox{ at } p_l + \Gamma, l = 1,2.
\end{displaymath}
As the poles have vanishing residues,
the function $\ \wp_1 - \wp_2\ $ integrates to a meromorphic function
$\ \tilde w \mbox{ on } \com\ $
with simple poles precisely at $\ p_1 + \Gamma, p_2 + \Gamma\ $
and with residues $\ (-1)^l \mbox{ at } p_l + \Gamma\ $.
As by (\ref{will-2.weier-residue})
\begin{displaymath}
	\int \limits_{[\xi,\xi + \omega_k]} (\wp_1 - \wp_2)(\zeta) \d \zeta = 0
\end{displaymath}
for $\ \omega_1 = 1, \omega_2 = \omega \mbox{ and appropriate } \xi \in \com\ $,
we see that $\ \tilde w\ $ is doubly periodic.
Comparing the residues of $\ w, \tilde w\ $,
we see that $\ w + \tilde w $ is doubly periodic meromorphic without poles,
hence constant. Therefore
\begin{displaymath}
	w' = -\tilde w' = -(\wp_1 - \wp_2).
\end{displaymath}
Putting
\begin{displaymath}
	\tilde f_3 := Re(-2 \tilde a_3 w),
\end{displaymath}
we calculate
\begin{displaymath}
	\partial_z \tilde f_3
	= -2 \partial_z Re(\tilde a_3 w)
	= -(\tilde a_3 w)'
	= \tilde a_3 (\wp_1 - \wp_2),
\end{displaymath}
which is (\ref{will-2.f3}).

If $\ \partial_z \tilde f(p) = 0 \mbox{ for some }
p \in \com - \{ p_1 + \Gamma , p_2 + \Gamma \}\ $,
then $\ \partial_z f_1(p), \partial_2 f_2(p) = 0\ $ and
\begin{displaymath}
	0 = \partial_z \tilde f_3 = \tilde a_3 (\wp_1 - \wp_2)(p) = 0.
\end{displaymath}
As $\ \tilde a_3 \neq 0\ $ by (\ref{will-2.tilde-not0}),
we get $\ (\wp_1 - \wp_2)(p) = 0\ $
and further by (\ref{will-2.fz}) and (\ref{will-2.ab-add}) that
\begin{displaymath}
	\partial_z f_j(p) = a_j \wp_1(p) + b_j \wp_2(p)
	= a_j (\wp_1 - \wp_2)(p) = 0
	\quad \mbox{for } j = 3,4.
\end{displaymath}
Together $\ \partial_z f(p) = 0\ $,
contrary to our assumption
that $\ f_0\ $ is an immersion on $\ \com\ $
and $\ f\ $ is an immersion on $\ \com - \{ p_1 + \Gamma, p_2 + \Gamma \}\ $.
Therefore $\ \partial_z \tilde f \neq 0 \mbox{ on } \com - \{ p_1 + \Gamma , p_2 + \Gamma \}\ $,
and
\begin{displaymath}
	\tilde f: \com - \{ p_1 + \Gamma, p_2 + \Gamma \} \rightarrow \rel^3
\end{displaymath}
is a conformal immersion.
Moreover $\ \tilde f_1 = f_1, \tilde f_2 = f_2
\mbox{ and } \tilde f_3 = Re(-2 \tilde a_3 w)\ $ are harmonic functions,
and therefore $\ \tilde f\ $ is a minimal conformal immersion.

%We proceed in proving that the inversion $\ \tilde f_0 := \tilde f / |\tilde f|^2:
%\coto{\omega} \rightarrow \rel^3\ $ is a smooth immersion,
%also at $\ p_1, p_2 \in \coto{\omega}\ $.
To analyse
%$\ \tilde f_0 \mbox{ and }
$\ \tilde f\ $ close to $\ p_1, p_2\ $,
we put $\ \hat f(z) := \tilde f(p_1 + (1/z))
\mbox{ for } |z| \gg 1\ $ and calculate
\begin{displaymath}
	\partial_z \hat f(z)
	= -\frac{1}{z^2} \partial_z \tilde f(p_1 + \frac{1}{z}) =
\end{displaymath}
\begin{displaymath}
	= -\frac{1}{z^2} \Big( \tilde a \wp(\frac{1}{z}) + \tilde b \wp(p_1 - p_2 + \frac{1}{z})
	+ \tilde d - (\tilde a + \tilde b) \wp(p_2-p_1) \Big) =
\end{displaymath}
\begin{displaymath}
	= -\frac{1}{z^2} \Big( \tilde a z^2 + O(1) \Big)
	= -\tilde a + O(\frac{1}{z^2}),
\end{displaymath}
hence
\begin{equation} \label{will-2.grad}
	\partial_z \hat f(z) \rightarrow \tilde a \neq 0
	\quad \mbox{for } z \rightarrow \infty.
\end{equation}
Therefore the gradient of $\ \hat f\ $ is bounded for large $\ z\ $,
and we can estimate
\begin{displaymath}
	|\hat f(z)| \leq C(1 + |z|).
\end{displaymath}
Putting $\ \hat f_\varrho(z) := \varrho \hat f(z/\varrho)\ $, we see
\begin{displaymath}
	|\hat f_\varrho(z)| \leq C (\varrho + |z|),
\end{displaymath}
hence $\ \hat f_\varrho\ $ is locally bounded in compact subsets of $\ \com - \{0\}\ $.
As further $\ \hat f_\varrho\ $ is harmonic,
we see for a subsequence that $\ \hat f_\varrho \rightarrow \hat f_0
\mbox{ smoothly in compact subsets of } \com - \{0\}\ $ and by above
\begin{displaymath}
	|\hat f_0(z)| \leq C |z|
\end{displaymath}
and
\begin{displaymath}
	\partial_z \hat f_0(z) \leftarrow \partial_z \hat f_\varrho(z) \rightarrow \tilde a,
\end{displaymath}
hence
\begin{displaymath}
	\hat f_0(z) = 2 Re(\tilde a) Re(z) - 2 Im(\tilde a) Im(z)
\end{displaymath}
and $\ \hat f_\varrho \rightarrow \hat f_0
\mbox{ smoothly in compact subsets of } \com - \{0\}\ $
for the whole family.
Clearly as $\ \langle \tilde a , \tilde a \rangle = 0, \tilde a \neq 0\ $,
we see that $\ Re(a) / |Re(a)| , Im(a) / |Im(a)|\ $ is an orthonormal basis
of a two dimensional subspace of $\ \rel^3\ $. 
%This yields for the image varifold
%$\ \tilde \mu^l := \Ht \llcorner \tilde f(B_\varrho(p_l)
%\mbox{ of } \tilde f \mbox{ for } \varrho > 0\ $ small
%that its blow down at infinity is given by
%\begin{equation} \label{will-2.blow-down}
%	T_\infty \tilde \mu^1 = span \{ Re(\tilde a) , Im(\tilde a) \},
%	\quad
%	T_\infty \tilde \mu^2 = span \{ Re(\tilde b) , Im(\tilde b) \},
%\end{equation}
%in particular $\ \theta^2(\tilde \mu^l,\infty) = 1\ $.
Therefore $\ \partial B_R(0)\ $ intersects $\ \tilde f(\coto{\omega} - \{p_1,p_2\})\ $
transversally in two closed curves whose preimages contract to $p_1$ resp. $p_2$ for large $\ R\ $,
and we get for the geodesic curvature
$\ \kappa_{\partial (\tilde f^{-1}(B_R(0)))}\ $ that
\begin{displaymath}
	\int \limits_{\tilde f^{-1}(\partial B_R(0)) }
	\kappa_{ \partial ( \tilde f^{-1}( B_R(0))) }  \d \sigma_{\tilde g}
	\rightarrow -4 \pi
	\quad \mbox{as } R \rightarrow \infty. 
\end{displaymath}
Here we let $\ \tilde g := \tilde f^* \geu \mbox{ on } \coto{\omega}\ $ be the pull-back metric and $\d \sigma_{\tilde g}$ is the induced line element. Moreover, we denote by $\ \mu_{\tilde g}\ $ the induced area measure and, since $\tilde f$ is minimal, the Gau{\ss}  curvature is given by $\ K_{\tilde g} = -|A_{\tilde f}|^2/2 \leq 0\ $, where $A_{\tilde f}$ is the second fundamental form of $\tilde f$.
Using the Gau\ss-Bonnet Theorem it follows that
\begin{displaymath}
	\int \limits_{\coto{\omega}-\{p_1,p_2\}} K_{\tilde g} \d \mu_{\tilde g}
	\leftarrow \int \limits_{\tilde f^{-1}(B_R(0))} K_{\tilde g} \d \mu_{\tilde g}=
\end{displaymath}
\begin{displaymath}
	= 2 \pi \chi(\tilde f^{-1}(B_R(0)))
	+\int \limits_{\tilde f^{-1}(\partial B_R(0)) }
	\kappa_{ \partial (\tilde f^{-1}( B_R(0))) }  \d \sigma_{\tilde g}	\rightarrow -8 \pi,
\end{displaymath}
where we used that $ \tilde f^{-1}(\partial B_R(0))$ consists of two closed curves which contract to $p_1$ resp. $p_2$ for $\ R \to \infty \ $, therefore these curves bound disks and hence $\tilde f^{-1}(B_R(0))\cong \coto{\omega} - \{p_1,p_2\}$ for $R$ large enough.
In particular, we obtain
\begin{displaymath}
	\int \limits_{\coto{\omega}-\{p_1,p_2\}} |A_{\tilde f}|^2 \d \mu_{\tilde g}=16\pi 
	< \infty.
\end{displaymath}
Then by \bcite{schoe.mini} Proposition 1 and Theorem 3,
we get that $\ \tilde f\ $ parametrizes a pair of planes or a catenoid.
But this is impossible, as $\ \tilde f $ is defined on $\ \coto{\omega} - \{p_1,p_2\}\ $,
and hence excludes the final case.
Therefore there is no immersion of a torus
with Willmore energy $\ 8 \pi\ $ and at least one double point.
\proof
In the next Proposition we show that the proof of the above theorem can be modified in order to classify all branched conformal immersions from a torus into $\rel^4$ with at least one branch point and Willmore energy $8\pi$.
\begin{proposition}\label{branch}
Any branched conformal immersion from
	$\coto{\omega}, \omega \in \modul$, into $\rel^4$
	with at least one branch point and with Willmore energy
	$8\pi$ is given, up to M\"obius transformations,
	by a branched double cover $\coto{\omega} \to S^2$,
	in particular, it has four branch points.
\end{proposition}
{\pr Proof:} 
We start by noting that the assumption on the Willmore energy implies that the multiplicity of every point $x\in f_0(\coto{\omega})$ can be at most two, see \cite{li.yau}. 
Next we assume that $p\in \coto{\omega}$ is a branch point of order two of $f_0$ and without loss of generality we let $f_0(p)=0$. As in the proof of the above theorem, we get that the inversion $f:=I(f_0)=f_0/|f_0|^2:\coto{\omega}-\{p\} \to \rel^4$ is a branched conformal immersion. Here we used that $f_0^{-1}(0)=\{p\}$ since otherwise the multiplicity of $0\in f_0(\coto{\omega})$ would be at least three. Additionally, it follows from Corollary 2 in \cite{nguy12} that $f$ is minimal. If we now choose local conformal coordinates around $p$ which map $p$ onto the origin, we have that $|f_0(z)| \asymp |z|^2$, $|\partial_z f_0(z)| \asymp |z|$. Hence we also get a local conformal parametrization of $f_0(B_1(0))$ and $f(B_1(0))$ around $p$ with
\begin{displaymath}
g(z):= (f^* \geu)(z) =(f_0^* I^* \geu)(z) = (f_0^* (|\cdot|^{-4} \geu))(z)=
\end{displaymath}
\begin{displaymath}
=|f_0(z)|^{-4} |\partial_z f_0(z)|^2 \geu \asymp |z|^{-6} \geu
\end{displaymath}
at least for $z$ close to $0$.
Since on the other hand
\[
g(z)=|\partial_z f(z)|^2 \geu  
\]
we conclude that 
\[
|\partial_z f| \asymp |z|^{-3}
\]
and hence $\partial f$ is a meromorphic $\com^4$-valued $1$-form on $\coto{\omega}$ with a pole of order three at $p$. Without loss of generality we assume from now on that $p=0$. Again, we consider $f_0$ and $f$ as branched conformal immersions which are doubly periodic with respect to $\Gamma$ on $\com - \Gamma$. 

Next we choose $a,b\in \com^4$ so that $\partial_z f -a\wp'-b\wp $ has a pole of order at most one in $\Gamma$. Note that $a\not= 0$.
Since there are no meromorphic doubly periodic functions with only one pole of order one, there exists $d\in \com^4$ so that
\[
\partial_z f = a\wp' +b\wp +d .
\]
Since $f$ is conformal away from its branch points, we conclude
\begin{displaymath}
	0 = \langle \partial_z f , \partial_z f \rangle =
\end{displaymath}
\begin{displaymath}
	= \langle a,a \rangle (\wp')^2 +  2\langle a, b \rangle \wp' \wp+ \langle b,b \rangle \wp^2 +
\end{displaymath}
\begin{displaymath}
	+ 2 \langle a,d \rangle \wp' + 2 \langle b,d \rangle \wp + \langle d,d \rangle.
\end{displaymath}
Comparing the order of the poles at the origin of the functions on the right hand side implies iteratively
\begin{equation}\label{coeffbranch}
0= \langle a,a \rangle = \langle a, b \rangle= \langle b,b \rangle= \langle a,d \rangle= \langle b,d\rangle= \langle d,d\rangle. 
\end{equation}
Arguing as in the previous theorem, as $a\not=0$, we can assume that $a=e_1-ie_2$ and therefore also $b_1-ib_2=0$ and $b_3-ib_4=0$ after a reflection. Moreover, $d_1-id_2=0$ and $(d_3+id_4)(d_3-id_4)=d_3^2+d_4^2=-d_1^2-d_2^2=0$. Lastly, we calculate
\[
0=\langle b,d\rangle =b_3(d_3-id_4).
\]
Now either $d_3-id_4=0$, from which we conclude 
\[
0=\partial_z f_1 -i\partial_z f_2 =\partial_z f_3-i\partial_z f_4
\]
and as in the previous theorem, this implies that the doubly periodic meromorphic functions $h_1=f_1+if_2$, $h_2=f_3+if_4$ have at most one double pole at $0$. In the second case we have $b_3=0$, which implies $b_4=0$ and $d_3+id_4=0$. Here we get
\[
0=\partial_z f_1 -i\partial_z f_2 =\partial_z f_3+i\partial_z f_4
\]
and this time we let $h_1=f_1+if_2$ and $h_2:=f_3-if_4$ and then we can repeat the above argument.

In both cases we conclude that there exist $\alpha_l,\beta_l\in \com$, $l=1,2$, with
\[
h_l=\alpha_l \wp +\beta_l
\]
and
\[
\partial_z f= \frac12 (h_1',-ih_1',h_2',-ih_2').
\]
This implies in particular that
\[
\partial_z f_3 = \partial_z Re(\alpha_2 \wp + \beta_2)
	= (1/2) (\alpha_2 \wp)' = (\alpha_2/2) \wp'.
\]
Since on the other hand we have $a_3=0$ this yields
\[
b_3 \wp +d_3 = 	(\alpha_2/2) \wp'
\]
and by comparing again the order of the poles of the involved functions we get
\[
b_3=b_4=d_3=d_4=\alpha_2=0.
\]
Hence we conclude that $f_3$ and $f_4$ are constant and $f_1=Re(\alpha_1 \wp+\beta_1)$, $f_2=Im(\alpha_1 \wp+\beta_1)$. Clearly $\alpha_1\not=0$ as $a\not= 0$ and  we can compose $f$ with a dilation and a translation in order to get that
\[
f:\coto{\omega} \to \rel^2 \times \{0\},\ \ \ f(z)=(\wp(z),0)
\]
and $f$ resp. $f_0$ is a branched double cover. Note that the degree of $\wp$ is two and the theorem of Riemann-Hurwitz, see \cite{jo.rie} Theorem 2.5.2, implies that $\wp$ has four branch points. 
\proof  

%%%%%

%%%%%
%%%%%

\setcounter{equation}{0}

\section{Infimal Willmore energy for fixed conformal class
in higher codimension} \label{infi}

In the previous section, we showed that there are no
immersed tori in $\ \rel^4\ $ with Willmore energy
equal to $\ 8 \pi\ $ and at least one double point. Additionally, we showed that the only branched conformal immersion with Willmore energy $8\pi$ is, up to M\"obius transformations, a branched double cover $T^2\to S^2$.
In this section, we show by perturbing
this branched immersion that
the infimal Willmore energy of smooth conformal immersions
is at most $\ 8 \pi\ $ in any conformal class of tori.
\begin{theorem} \label{infi.tori}

%%%%%
For any conformal class $\ \omega \in \modul\ $, we have
\begin{equation} \label{infi.tori.esti}
	\miniz{4}(\omega) \leq 8 \pi,
\end{equation}
in particular $\ \minin \mbox{ is continuous for } n \geq 4\ $.%%%%%
\end{theorem}
{\pr Proof:} \\
For $\ \sigma \mbox{ close to } \omega\ $ the Weierstrass $\ \wp-$function is given by
\begin{displaymath}
	\wp_\sigma(z) := \frac{1}{z^2} + \sum \limits_{\gamma \in \ganz + \sigma \ganz - \{0\}}
	\Big( \frac{1}{(z-\gamma)^2} - \frac{1}{\gamma^2} \Big),
\end{displaymath}
see the proof of Proposition \ref{nonexistence} and \bcite{ahl} \S 7.3.1.
This is a doubly-periodic, meromorphic function
with one pole of order two at the origin
apart from congruences with respect to the lattice $\ \ganz + \sigma \ganz\ $.
It can be considered as a holomorphic mapping
$\ \coto{\sigma} = \com / \ganz + \sigma \ganz \rightarrow S^2\ $,
which is of degree two and has four branch points
according to the theorem of Riemann-Hurwitz, see \bcite{jo.rie} Theorem 2.5.2.
To work with simple poles as in \S \ref{will},
we choose $\ \alpha \in \com\ $ which has two preimages
with non-vanishing derivative of $\ \wp_\omega\ $
and see that $\ (\wp_\sigma - \alpha)^{-1}\ $
is doubly-periodic and meromorphic
with two simple poles apart from congruences for $\ \sigma \approx \omega\ $.
We transform this on the fixed reference torus $\ T^2 := \com / \ganz + i \ganz\ $
by the linear map $\ A_\sigma: \com \diff \com, A_\sigma := (e_1 , \sigma)\ $,
which maps $\ \ganz + i \ganz \mbox{ onto } \ganz + \sigma \ganz\ $,
hence factors to a diffeomorphism $\ T^2 \diff \coto{\sigma}\ $,
and put
\begin{displaymath}
	f_\sigma := (\wp_\sigma - \alpha)^{-1} \circ A_\sigma \circ (\com \rightarrow T^2)^{-1}:
	T^2 \rightarrow \com \cup \{ \infty \}
\end{displaymath}
for the canonical quotient map $\ \com \rightarrow T^2\ $.
As the poles are simple now,
we know from Proposition \ref{will.prop}
that for any M\"obius transformation $\ \Phi: \rel^4 \cup \{ \infty \} \diff S^4\ $
the map $\ \Phi \circ (f_\sigma,0): T^2 \rightarrow S^4\ $ is smooth
and an immersion outside the branch points of $\ f_\sigma\ $.
We note that $(\wp_\sigma(z)-\alpha)^{-1}$ is a continuous function which is analytic in the variables $z$, $\sigma$ separately away from its poles. Therefore, it follows from Osgood's Lemma (see e.g. \cite{GuRo}), that $\ f_\sigma\ $ depends analytically
on $\ (z,\sigma)\ $ outside the poles.
Clearly, $\ f_\sigma\ $ induces the conformal class $\ \sigma\ $,
provided that $\ \sigma \in \modul\ $.
To avoid problems when $\ \omega \in \partial \modul\ $,
we work instead in the Teichm\"uller space $\ \teich = [Im > 0] \subseteq \com\ $,
then clearly $\ f_\sigma\ $ induces the Teichm\"uller class $\ \sigma\ $.
For the Teichm\"uller space,
we refer to \bcite{fisch.trom.conf} and \bcite{trom.teich}.

Let $\ p_0, p_1 \in T^2\ $ be the two poles
and $\ b_1, \ldots b_4 \in T^2\ $ be the branch points of $\ f_\omega\ $,
which are distinct, as the poles are simple.
We fix $\ \delta > 0\ $ such that the balls $\ B_{3 \delta}(p_i), B_{3 \delta}(b_j)
\mbox{ for } i = 0,1, j = 1,2,3,4,\ $ are pairwise disjoint
congruent to $\ \ganz + i \ganz\ $
and are evenly covered by $\ \com \rightarrow T^2\ $.
Next for $\ \sigma \in \overline{B_\varrho(\omega)}
\mbox{ with } \varrho > 0\ $ small,
$\ f_\sigma\ $ has exactly one pole respectively one branch point
in $\ B_\delta(p_i), B_\delta(b_j)\ $.
As in \S \ref{will}, we add a function $\ \varphi_\sigma\ $
in order that $\ \Phi \circ (f_\sigma , \varepsilon \varphi_\sigma)\ $
is an immersion also at $\ b_j \mbox{ for } \varepsilon > 0\ $ and for every M\"obius transformation $\ \Phi:\rel^4 \cup \{ \infty \} \diff S^4 \ $.
This is done by choosing $\ \varphi_\sigma\ $ to be holomorphic in $\ B_{2 \delta}(b_j)\ $,
more precisely we consider $\ z = (x,y) \mapsto A_\sigma z\ $
which obviously is holomorphic in the complex structure induced
by $\ A_\sigma: \com \diff \com\ $.
Then we put for some cut-off function $\ \eta \in C^\infty_0(B_{3 \delta}(0)),
0 \leq \eta \leq 1, \eta \equiv 1 \mbox{ in } B_{2 \delta}(0)\ $, that
\begin{displaymath}
	\varphi_\sigma(p) := \sum \limits_{j=1}^4
	\eta(p - b_j) A_\sigma p
	\quad \mbox{for } p \in T^2,
\end{displaymath}
where we take any branch of $\ A_\sigma \mbox{ locally in } B_{3 \delta}(b_j)\ $.
We calculate in $\ B_{2 \delta}(b_j)\ $ by holomorphy
\begin{displaymath}
	(f_\sigma , \varepsilon \varphi_\sigma)^* \geu
	= ( (\wp_\sigma - \alpha)^{-1} \circ A_\sigma , \varepsilon A_\sigma)^* \geu
	= A_\sigma^* ( |(\wp_\sigma - \alpha)^{-1})'|^2 + \varepsilon^2) \geu
\end{displaymath}
and see that $\ (f_\sigma , \varepsilon \varphi_\sigma)\ $
is a smooth immersion in $\ B_{2 \delta}(b_j) \mbox{ for } \varepsilon > 0\ $
which is conformal to $\ A_\sigma^* \geu\ $. Similarly, we have outside of $\ \cup_{j=1}^4 B_{3 \delta}(b_j) \cup  \{p_1,p_2\} \ $ 
\begin{displaymath}
	(f_\sigma , \varepsilon \varphi_\sigma)^* \geu
	= ( (\wp_\sigma - \alpha)^{-1} \circ A_\sigma , 0)^* \geu
	= A_\sigma^*  |(\wp_\sigma - \alpha)^{-1})'|^2  \geu.
\end{displaymath}
As $\ \Phi \circ f_\sigma\ $ is an immersion in $\ T^2 - \cup_{j=1}^4 B_\delta(b_j)\ $
and $\ supp\ \varphi_\sigma \subseteq \cup_{j=1}^4 B_{3 \delta}(b_j)\ $,
we see that  Proposition \ref{will.prop} implies that $\ \Phi \circ (f_\sigma , \varepsilon \varphi_\sigma)\ $ is an immersion on $\ T^2 \mbox{ for } \varepsilon > 0\ $, $\ \sigma \in \overline{B_\varrho(\omega)}
\mbox{ for } \varrho > 0\ $ small
and for any M\"obius transformation $\ \Phi: \rel^4 \cup \{ \infty \} \diff S^4\ $.
Moreover the pull-back metric of $\ \Phi \circ (f_\sigma , \varepsilon \varphi_\sigma)\ $
is conformal to $\ A_\sigma^* \geu\ $ outside of
\begin{displaymath}
	\cup_{j=1}^4 [d \eta(\cdot - b_j) \neq 0]
	\subseteq \bigcup_{j=1}^4 \Big( B_{3 \delta}(b_j) - B_{2 \delta}(b_j) \Big).
\end{displaymath}
Next we introduce the projection
$\ \tp: \metric = \{ \mbox{ smooth metrics on } T^2 \} \rightarrow \teich\ $
in the Teichm\"uller space
and examine the Teichm\"uller class induced
by $\ \Phi \circ (f_\sigma , \varepsilon \varphi_\sigma)\ $ by putting for every $\varepsilon>0$
\begin{displaymath}
	\tau(\sigma,\varepsilon) :=  \pi \Big( (\Phi \circ (f_\sigma , \varepsilon \varphi_\sigma))^* \geu \Big).
%	= \pi \Big( (f_\sigma , \varepsilon \varphi_\sigma)^* \geu \Big)
\end{displaymath}
The projection $\ \tp\ $ is smooth in the $\ W^{s,2}-\mbox{topology of } \metric\ $,
see \bcite{fisch.trom.conf}
Theorem 1.5, 2.2, 7.7, 7.8, \S 8,9 for $\ s > 2\ $
and \bcite{trom.teich}
Theorem 1.3.2, 1.6.2, Corollary 1.3.3,
Remark 2.5.3 (1) for $\ s > 3\ $. Therefore it follows from the above smoothness considerations for $f_\sigma$ that $\tau$ is continuous in $\mathcal T \times \rel^+$ and
we claim that $\ \lim_{n\to \infty} \tau(\sigma_n,\varepsilon_n) = \sigma \ $ as $\ (\sigma_n,\varepsilon_n)\to (\sigma,0) \ $. In order to see this we let $\ \eta \ $ be as above and we define the function $\ \lambda_{\sigma, \varepsilon} : \com \to \rel^+ \ $ by 
\begin{displaymath}
\lambda_{\sigma,\varepsilon}(z)= 1+ \sum_{j=1}^4 \eta(z-b_j) \left( \frac{1}{|(\wp_\sigma - \alpha)^{-1}(A_\sigma(z)))'|^2 + \varepsilon^2} -1\right)
\end{displaymath}
\begin{displaymath}
+\sum_{i=1}^2 \eta(z-p_i) \left( \frac{1}{|(\wp_\sigma - \alpha)^{-1}(A_\sigma(z)))'|^2} -1\right).
\end{displaymath}
Using this definition and the above expressions for the metric $\ (f_\sigma , \varepsilon \varphi_\sigma)^* \geu \ $ near the branch points resp. poles, it follows that the new metric
\begin{displaymath}
 \lambda_{\sigma,\varepsilon} (f_\sigma , \varepsilon \varphi_\sigma)^* \geu
\end{displaymath}
can be smoothly extended into the poles and is independent of $\varepsilon$ outside of $\ \cup_{j=1}^4 B_{3\delta} (b_j) - B_{2\delta}(b_j)\ $. Moreover, it depends smoothly on $\ (\sigma,\varepsilon) \ $ and is a smooth nondegenerate metric also for $\ \varepsilon =0 \ $. Therefore, the new metric converges smoothly and since 
\begin{displaymath}
\tau(\sigma,\varepsilon)
= \pi( \lambda_{\sigma,\varepsilon}
(f_\sigma,\varepsilon \varphi_\sigma)^* \geu )
\end{displaymath}
and 
\begin{displaymath}
\pi ( \lambda_{\sigma,0}
(f_\sigma,0)^* \geu ) =\sigma
\end{displaymath}
we finish the proof of the claim.
As by the above $\ \lim_{\varepsilon \to 0} \tau(\sigma,\varepsilon) = \sigma \ $,
we know for $\ \varepsilon > 0\ $ small that the mapping degree is
\begin{displaymath}
	deg( \tau(\cdot,\varepsilon) , \overline{B_\varrho(\omega)} , \omega ) = 1,
\end{displaymath}
hence there exists $\ \sigma_\varepsilon \approx \omega
\mbox{ with } \tau(\sigma_\varepsilon,\varepsilon) = \omega\ $
and therefore, for every $\ \varepsilon > 0\ $ small, we have
\begin{displaymath}
	\miniz{4}(\omega)
	\leq \W( \Phi \circ (f_{\sigma_\varepsilon} , \varepsilon \varphi_{\sigma_\varepsilon}) ).
\end{displaymath}
As in the previous proposition, it follows from Gauss-Bonnet that
\begin{displaymath}
\W ( \Phi \circ (f_{\sigma_\varepsilon} , \varepsilon \varphi_{\sigma_\varepsilon}) )=8\pi +\W (  (f_{\sigma_\varepsilon} , \varepsilon \varphi_{\sigma_\varepsilon}) ).
\end{displaymath}
Since $\ (f_{\sigma_\varepsilon} , \varepsilon \varphi_{\sigma_\varepsilon})\ $ is holomorphic and hence minimal outside of $\ \bigcup_{j=1}^4 \Big( B_{3 \delta}(b_j) - B_{2 \delta}(b_j)\Big) \ $ we get
\begin{displaymath}
\W ( (f_{\sigma_\varepsilon} , \varepsilon \varphi_{\sigma_\varepsilon}))
= \W (  (f_{\sigma_\varepsilon} , \varepsilon \varphi_{\sigma_\varepsilon}) |_{\bigcup_{j=1}^4 (B_{3 \delta}(b_j) - B_{2 \delta}(b_j) )}).
\end{displaymath}
Moreover, $\ (f_{\sigma_\varepsilon} , \varepsilon \varphi_{\sigma_\varepsilon}) \to (f_\sigma,0) \ $ smoothly on $\ \bigcup_{j=1}^4 \Big( B_{3 \delta}(b_j) - B_{2 \delta}(b_j)\Big) \ $, which implies
\begin{displaymath}
\W (  (f_{\sigma_\varepsilon} , \varepsilon \varphi_{\sigma_\varepsilon}) |_{\bigcup_{j=1}^4 (B_{3 \delta}(b_j) - B_{2 \delta}(b_j) )} )
\rightarrow \W (  (f_\sigma,0))
|_{\bigcup_{j=1}^4 (B_{3 \delta}(b_j) - B_{2 \delta}(b_j) )})
= 0
\end{displaymath}
and therefore
\begin{displaymath}
	\miniz{4}(\omega)
	\leq	\lim \limits_{\varepsilon \rightarrow 0}
	\W( \Phi \circ (f_{\sigma_\varepsilon} , \varepsilon \varphi_{\sigma_\varepsilon}) )
	= 8 \pi,
\end{displaymath}
thereby establishing (\ref{infi.tori.esti}).

The continuity of $\ \minin \mbox{ for } n \geq 4\ $
follows directly from (\ref{infi.tori.esti})
and \bcite{schae.comp-will12} Proposition 4.1.
\proof
{\large \bf Remark:} \\
The above construction is not possible in codimension one.
Indeed, let $\ f_k: \coto{\omega} \rightarrow \rel^3,
\omega \in \modul\ $, be a sequence of conformal immersions with
\begin{displaymath}
	\limsup \limits_{k \rightarrow \infty} \W(f_k) \leq 8 \pi
\end{displaymath}
and
\begin{displaymath}
	f_k \rightarrow f_0 \mbox{ weakly in } W^{2,2}_{loc}(\coto{\omega} - S,\rel^3)
\end{displaymath}
for some finite set $\ S \subseteq \coto{\omega}\ $
and where $\ f_0: \coto{\omega} \rightarrow \rel^3\ $
is a branched conformal $\ W^{2,2}-$immersion with square integrable second fundamental form.
We know that this is true by \cite{kuw.li.mini} Proposition 4.1
for a subsequence $\ f_k \mbox{ replaced by } \Phi_k \circ f_k\ $
for appropriate M\"obius transformations $\ \Phi_k \mbox{ of } \rel^3\ $.

Now we assume that $\ f_0\ $ has at least one branch point $\ p \in \coto{\omega}\ $.
We have $\ \W(f_0) \leq 8 \pi\ $ by lower semicontinuity,
and assuming $\ f_0(p) = 0\ $ after translation,
we get as in the proof of Proposition \ref{branch}
with \cite{li.yau} and \cite{nguy12} Corollary 2
that the inversion $\ f := I(f_0) = f_0/|f_0|^2: \coto{\omega} - \{p\} \to \rel^3\ $
is a branched conformal minimal immersion.
Then outside the finitely many branch points of $\ f_0 \mbox{ respectively of } f\ $,
hence almost everywhere on $\ \coto{\omega} - \{p\}\ $,
we get $\ \Delta f = 0\ $, in particular $\ f \mbox{ is smooth on } \coto{\omega} - \{p\}\ $.
Moreover $\ |f_0(p+z)| \asymp |z|^2 \mbox{ and } |\nabla f_0(p+z)| \asymp |z|\ $
by \cite{kuw.li.mini} Theorem 3.1 and its proof,
as the branch point $\ p \mbox{ of } f_0\ $ has order two,
hence as in the proof of Proposition \ref{branch}
\begin{displaymath}
	|\partial_z f(p + z)| \asymp |z|^{-3}
\end{displaymath}
and $ \partial f\ $ is a meromorphic $\ \com^3$-valued $1$-form
on $\ \coto{\omega}\ $ with a pole of order three at $\ p\ $.
Then the proof of Proposition \ref{branch} proceeds,
and $\ f_0\ $ is up to M\"obius transformations
a branched double cover $\ \coto{\omega} \to S^2\ $,
in particular it has four branch points.
By the proof of \cite{kuw.li.mini} Proposition 4.1,
these branch points are contained in the exceptional set $\ S\ $,
hence $\ \#(S) \geq 4\ $, and
\begin{displaymath}
	\int \limits_{\coto{\omega}} |A_{f_0}|^2 \d \mu_{f_0}
	+ \#(S) \gamma_3
	\leq \limsup \limits_{k \rightarrow \infty} \int \limits_{\coto{\omega}}
	|A_{f_k}|^2 \d \mu_{f_k},
\end{displaymath}
where $\ \gamma_3 = 8 \pi\ $ in codimension one,
see \cite{kuw.li.mini} Corollary 2.4 and Proposition 4.1
and \cite{schae.comp-will12} Proposition 5.1.
Since by the Gau\ss-Bonnet theorem and the Gau\ss\ equations
\begin{displaymath}
	\limsup \limits_{k \rightarrow \infty}
	\int \limits_{\coto{\omega}} |A_{f_k}|^2 \d \mu_{f_k}
	= \limsup \limits_{k \rightarrow \infty}
	\Big( 4 \W(f_k) - 2 \int \limits_{\coto{\omega}} K_{f_k} \d \mu_{f_k} \Big)
	\leq 32 \pi
\end{displaymath}
and $\ \#(S) \geq 4\ $, we conclude that $\ A_{f_0} \equiv 0\ $,
which is not true, as $\ f_0: \coto{\omega} \rightarrow S^2\ $ is a branched double cover.

Therefore $\ f_0\ $ is an unbranched conformal $\ W^{2,2}-$immersion.
Moreover by \cite{kuw.li.mini} Theorem 3.1,
we see that $\ f_0\ $ is uniformly conformal
in the sense that $\ f_0^* \geu = e^{2u} \geu
\mbox{ with } u \in L^\infty(\coto{\omega})\ $.
If further $\ \miniz{3}(\omega) = 8 \pi\ $,
that is $\ f_k\ $ is a minimizing sequence,
then by \cite{kuw.schae.will7} Theorem 7.4,
we get that $\ f_0\ $ is a smooth conformally constrained minimizer
of the Willmore energy on $\ \coto{\omega} \mbox{ in } \rel^3\ $.

\defin
\\ \ \\
Combining with \cite{kuw.li.mini} Corollary 4.1 or \cite{rivi.mini} Theorem 1.17,
we obtain the following corollary
which improves the existence result for conformally constrained Willmore minimizers
of Kuwert-Li \cite{kuw.li.mini} and Rivi\`ere \cite{rivi.mini} in the case of tori in $\rel^3$.

\begin{corollaryth}
For every conformal class $\omega \in \modul$ with $\ \miniz{3}(\omega) \leq 8 \pi\ $,
there exists a smooth conformal immersion $f:\coto{\omega} \to \rel^3$
which minimizes the Willmore energy
in the set of all conformal immersions
on $\ \coto{\omega} \to \rel^3\ $.
\end{corollaryth}

%%%%%

%\newpage

%\input{app}
%%%%%

%%%%%

\end{document}